\documentclass[a4paper,12pt]{article}
\usepackage{amsfonts,amsthm,amsmath,amssymb,graphicx}

\title{Dimension estimates for invariant measures of contracting-on-average iterated function systems}

\date{}

\author{Micha\l \ Rams \thanks{supported by Polish KBN Grant No 2P0 3A 034 25}\\
\normalsize Institute of Mathematics, Polish Academy of Sciences\\
\normalsize ul. \'Sniadeckich 8, 00-950 Warszawa, Poland\\
\normalsize e-mail rams@impan.gov.pl}

\theoremstyle{plain}
\newtheorem{lem}{Lemma}[section]
\newtheorem{prop}[lem]{Proposition}
\newtheorem{thm}[lem]{Theorem}
\newtheorem{cor}[lem]{Corollary}

\theoremstyle{definition}

\newtheorem{ex}[lem]{Example}

\theoremstyle{remark}

\numberwithin{equation}{section}

\DeclareMathOperator{\diam}{diam}
\DeclareMathOperator{\vol}{vol}
\DeclareMathOperator{\dist}{dist}
\DeclareMathOperator{\var}{var}

\DeclareMathOperator{\id}{Id}

\newcommand{\R}{\mathbb R}
\newcommand{\Z}{\mathbb Z}
\newcommand{\N}{\mathbb N}

\renewcommand{\epsilon}{\varepsilon}

\begin{document}

\maketitle

\begin{abstract}
We estimate from above and below the dimension of invariant measure for 
contracting-on-average iterated function systems in $\R^d$.
\end{abstract}

\def\thefootnote{}
\footnote{1991 {\it Mathematics
Subject Classification}: Primary 28A78,34D45}
\def\thefootnote{\arabic{footnote}}
\setcounter{footnote}{0}

\section{Introduction}

The goal of this paper is to give dimension estimates for invariant 
measures of contracting-on-average self-conformal iterated function
systems with point-depending probabilities.

The motivation for this work comes from the recent results of Szarek and Myjak \cite{S}, \cite{MS} and also Nicol, Sidorov and Broomhead \cite{NSB} and Fan, Simon and Toth \cite{FST} who gave bounds for the Hausdorff dimension of the invariant measure for some contracting-on-average iterated function systems.

Lately, Johansson and \"Oberg computed the precise value of the Hausdorff dimension for bounded (not necessarily contracting) iterated function systems on the line satisfying the open set condition, \cite{JO}.
Similar questions come from other branches of mathematics, compare \cite{Wa}.

There are two kinds of results on the dimension of an invariant measure for (contracting) iterating function systems.
The upper bound is easy to obtain, depending on assumptions it takes form of Moran formula (for similitudes) or Bowen (pressure) formula (for smooth conformal maps).
For systems with probabilities we get Eckmann-Ruelle formula, i.e. ratio of entropy to the Lyapunov exponent.
To get the lower bound for an iterated function system one typically needs a version of the open set condition (like in \cite{LN}).

We are working here with systems that aren't everywhere contracting.
We need to assume some contractiveness-on-average condition to guarantee existence and uniqueness of the invariant measure, as was proved in \cite{BDEG}.

We prove the Eckmann-Ruelle formula for the upper bound on the Hausdorff dimension of the invariant measure.
It is the same as the bound in \cite{JO} and \cite{FST}, but we work on different class of systems, in particular our method works for measures with unbounded supports.
As for the lower bound, we can only obtain it under some (bit stronger) version of OSC.

The rest of the paper is divided as follows.
In the second section we state some known results on the iterated function systems and formulate our theorem.
In the third section we prove some properties of contracting-on-average iterated function systems that we will need in the future. 
As the results of this section might be of independent interest, they are formulated in more general situation, i.e. for iterated function system acting on a Polish space and satisfying Dini condition.
In the fourth section we finish our preparation by formulating some geometric lemmas about typical behaviour of orbits of IFS.
Finally, in the fifth section we present the proof of Theorem \ref{thm:main}.

\section{Introduction to iterated function systems}

For informations on dimension theory and properties of contracting iterated function systems we refer the reader to standard textbooks like \cite{F}.
The existence and most important properties of invariant measures for contracting-on-average iterated function systems are to be found in \cite{BDEG} and \cite{E}.

Consider an iterated function system consisting of $k$ $C^{1+\alpha}$ maps $h_1,\ldots,h_k:\R^d\rightarrow \R^d$ and a probability vector field $\vec{p}(x)=(p_1(x),\ldots,p_k(x))$ defined on $\R^d$.
For a map $h$ we denote $h'(x)=||Dh(x)||$.
We assume $h_i'$ are uniformly bounded away from 0 and $\infty$.
We also assume that $p_i$ are $C^\alpha$ and uniformly bounded away from 0.
Hence, $\log p_i$ and $\log h_i'$ are $C^{\alpha}$.

The system acts as follows: on every iteration we randomly choose $i\in\{1,\ldots,k\}$ according to the vector $\vec{p}(x)$ (depending on the point $x\in \R^d$ we are in) and then apply the map $h_i$.
If all the maps $h_i$ are contracting, there exists a unique probability measure solving equation
\begin{equation} \label{eqn:eqn}
\nu = F(\nu):=\sum_i h_{i*} (p_i \cdot \nu),
\end{equation}
it is called {\it invariant measure} for the iterated function system.
Moreover, the invariant measure satisfies
\begin{equation} \label{eqn:ifs}
\overline{\dim}_H (\nu) \leq \frac {\int \sum_i p_i(x) \log p_i(x) d\nu(x)} {\int \sum_i p_i(x) \log |h_i'(x)| d\nu(x)}
\end{equation}
where the upper Hausdorff dimension of the measure is defined as the infimum of the Hausdorff dimension over sets of full measure.
Let us remind here that the lower Hausdorff dimension $\underline{\dim}_H$ is defined as the infimum of the Hausdorff dimensions of sets of positive measure and if both are equal, the common value is called Hausdorff dimension of the measure and denoted $\dim_H$.

Denote $\Sigma=\{1,\ldots,k\}^\N$.
Given a sequence $\omega\in\Sigma$, we will denote by $\omega_n$ the
$n$-th element and by $\omega^n$ the first $n$ elements of $\omega$.
We will also denote $\omega^{-n}=\omega_n \omega_{n-1} \ldots \omega_1$.
The {\it $n$-th level cylinder} $C_{\omega^n}$ is the set of all sequences from $\Sigma$ that begin with $\omega^n$.
Let $\sigma, \sigma_i$ be the left shift map and the $k$ right shift maps on $\Sigma$, i.e. 
\[\sigma \omega = \omega_2 \omega_3 \ldots\]
\[\sigma_i \omega = i \omega_1 \omega_2 \ldots\]

Define $h_{\omega^n}=h_{\omega_n}\circ\ldots\circ h_{\omega_1}$.
Let 
\[p_{\omega^n}(x)=\prod_{i=1}^n p_{\omega_i}(h_{\omega^{i-1}}(x))\]
and let $p_x$ be the unique probability measure on $\Sigma$ satisfying 
\[p_x(C_{\omega^n})=p_{\omega^n}(x)\]
The limit
\begin{equation} \label{eqn:pidef}
\pi(\omega)=\lim_{n\rightarrow\infty} h_{\omega^{-n}}(x_0)
\end{equation}
exists for all $\omega\in\Sigma$ and does not depend on $x_0\in \R^d$.
The invariant measure $\nu$ can be written as
\begin{equation} \label{eqn:mudef}
\nu=\pi_*(\mu_-)
\end{equation}
where $\mu_-$ is some $\sigma$-invariant probability measure on $\Sigma$.

In this paper we assume only that the IFS is {\it contracting-on-average}, i.e. 
\begin{equation} \label{eqn:coa}
\sum_1^k p_j(x) \log \frac {d(h_j(x),h_j(y))} {d(x,y)}<-b<0
\end{equation}
uniformly in $x$ and $y$.
The existence and uniqueness of the invariant measure $\nu$ for such systems was established in \cite{BDEG}.
In addition, it was proved that the measure is attractive, i.e. that $F^n(\tilde{\nu})\rightarrow \nu$ for all probability measures $\tilde{\nu}$.

To get a reasonable lower bound on the Hausdorff dimension of the invariant measure we will need stronger geometric assumptions.
Namely, we will assume $h_i$ are conformal maps.
The conformality assumption is void for $d=1$ but very strong in higher dimension: it implies analyticity in dimension 2 and only M\"obius transformations are conformal in $\R^d$, $d\geq 3$ (Liouville Theorem).

In conformal case, if IFS is contracting and if in addition there exists an open set $U$ such that $h_i(U)\subset U$ and $h_i(U)\cap h_j(U)=\emptyset$ for $i\neq j$ (it is called {\it open set condition}, OSC), we have equality in \eqref{eqn:ifs} (\cite{F}).

Our main result is that the results on the Hausdorff dimension of invariant measures for contracting-on-average iterated function systems are almost the same as for contracting ones.
The one difference is that instead of OSC we have to use a bit stronger geometrical conditions.

If the IFS satisfies OSC for some open set $U$ and if (for some $R_1$)
\begin{equation} \label{eqn:str}
\dist(h_i(U),h_j(U))\geq R_1>0,
\end{equation}
the system satisfies {\it strong open set condition} (SOSC).
If (in addition to OSC) $U$ satisfies (for some $R_2$ and $R_3$)
\begin{equation} \label{eqn:reg}
\forall_{r<R_2} \forall_{x\in U} \vol(B_r(x)\cap U)\geq R_3 r^d,
\end{equation}
the system satisfies {\it regular open set condition} (ROSC).

Denote the right hand side of \eqref{eqn:ifs} by $s$.

\begin{thm} \label{thm:main}
Let $(\R^d, h_i, p_i)$ be an iterated function system.
Assume $h_i$ are $C^\alpha$ and with bounded $\log h_i'$, $p_i$ are $C^\alpha$ and bounded away from 0.
Assume that \eqref{eqn:coa} is satisfied. 
Denote the invariant measure by $\nu$.
Then 
\[\overline{\dim}_H(\nu)\leq s\]

If in addition all $h_i$ are conformal and either SOSC or ROSC is satisfied, 
\[\dim_H(\nu)=s\]
\end{thm}

\section{Coding map and ergodicity}

Our main goal in this section is to prove the following proposition:

\begin{prop} \label{prop:pidef}
Let $X$ be a Polish space and fix $x_0\in X$.
Let $(X, h_i, p_i)$ be an iterated function system satisfying \eqref{eqn:coa}, with $p_i$ satisfying Dini condition and bounded away from 0 and with $h_i$ Lipschitz.
Let $\nu$ be the invariant measure of this system.
Then there exists a $\sigma$-invariant probabilistic measure $\mu_-$ on $\Sigma$ such that \eqref{eqn:pidef} is well defined $\mu_-$-almost everywhere and \eqref{eqn:mudef} holds. 
\end{prop}

Similar result was proved by Werner \cite{We}, but he used a stronger assumption instead of \eqref{eqn:coa}.

We start with some notations.
Let
\[\mu_+(d\omega)=\int_X p_x(d\omega) d\nu(x)\]
It is a probability $\sigma$-invariant measure on $\Sigma$.
Denote $\tilde{\Sigma}=\{1,\ldots,k\}^\Z$ and let $\tilde{\sigma}$ be the left shift on $\tilde{\Sigma}$, i.e. $(\tilde{\sigma}\omega)_i=\omega_{i+1}$.
The space $(\tilde{\Sigma},\tilde{\sigma})$ has a natural equivalence with $(\Sigma\times\Sigma,S)$, where
\[S(\tau,\omega)=(\sigma_{\omega_1}\tau, \sigma\omega)\]
On $\tilde{\Sigma}$ there exists a unique $\tilde{\sigma}$-invariant measure $\mu$ such that $\mu_+$ is the marginal measure of $\mu$ with respect to the positive coordinates:
\[\mu_+(\cdot) = \int \mu(d\tau, \cdot)\]
The measure $\mu$ is called {\it natural extension} of $\mu_+$.

We define the measure $\mu_-$ as the marginal measure of $\mu$ with respect to the nonpositive coordinates, i.e.
\[\mu_-(\cdot) = \int \mu(\cdot, d\omega)\]
The explicit formula for $\mu_-$ can be written as
\begin{equation} \label{eqn:mumi}
\mu_-(C_{\omega^n})=\mu_+(C_{\omega^{-n}})
\end{equation}

Let $\pi$ be given by \eqref{eqn:pidef} (wherever it is defined).
In order to prove \eqref{eqn:mudef} we need only to prove that $\pi$ is well defined for $\mu_-$-almost all $\omega\in\Sigma$.
Indeed, assume it is so.
The measure $\pi_*(\mu_-)$ is then a probability measure on $X$. 
Moreover, 
\[h_i(\pi(\omega))=\pi(\sigma_i \omega)\]
hence $\pi_*(\mu_-)$ is invariant under the action of the IFS.
By uniqueness of the invariant measure, it must be equal to $\nu$ and \eqref{eqn:mudef} holds.

We will need some estimations on the ergodic sums of integrable functions over typical trajectories of the IFS.
The following lemma is a weak large deviation estimation:

\begin{lem} \label{lem:clt}
Choose arbitrary bounded $A\subset X$.
Then there exists $L(A)$ such that for all $x\in X$, $n>0$ and $K>0$
\[p_x(\omega\in\Sigma; \diam(h_{\omega^n}(A))\leq K \exp(-nb))> 1-\frac {L(A)(1+\log^+(d(x,x_0)))n} {(\log K)^2} \]
\begin{proof}
For $j\in \{1,\ldots,k\}$ and $y\in X$ let 
\[\tilde{H}(y,j)= \sup_{z\in X} \log \frac {d(h_j(y),h_j(z))} {d(y,z)}\]
and
\[H(y,j)= \tilde{H}(y,j) -b - \sum_{i=1}^k p_i(y) \tilde{H}(y,i)\]
By \eqref{eqn:coa}, $H(y,j)> \tilde{H}(y,j)$ for all $y$ and $j$.
In addition, as $\tilde{H}$ is uniformly bounded, so is $H$.

Define a sequence of random variables of fixed expected value $b$:
\[H_n(\omega)=H(h_{\omega^{n-1}}(x),\omega_n)\]
where $\omega\in (\Sigma, p_x)$.
By boundedness of $H_n$, we have a bound on its variation as well:
\[\var H_n < v\]

The expected value of $H_n$ is the same, independently of $\omega^{n-1}$.
It implies that the random variables $H_n$ and $\sum_{j=1}^{n-1} H_j$ are uncorrelated (but not necessarily independent).
Hence
\[\var (\sum_{j=1}^n H_j) < vn\]
and by Tchebysheff inequality
\begin{equation} \label{eqn:pxo}
p_x(\omega; \sum_{j=1}^n H_j(\omega) > -nb + \log K) < \frac v {(\log K)^2}
\end{equation}

Now, we can estimate
\[\diam(h_{\omega^n}(A))\leq 2 \sup_{y\in A} d(h_{\omega^n}(x),h_{\omega^n}(y))\leq 2 \exp (\sum_{j=1}^n H_j(\omega)) (d(x,x_0)+d(x_0,A)+\diam A)\]
(here we use the inequality $\tilde{H}<H$ and \eqref{eqn:coa}).
The assertion follows from \eqref{eqn:pxo}.
\end{proof}
\end{lem}

\begin{lem} \label{lem:log}
\[\int_X \log^+(d(x,x_0)) d\nu(x)\leq M(x_0)<\infty\]
\begin{proof}
Let $\tilde{\nu}$ be any probability measure on $X$ with bounded support.
Let 
\[M_0=\max_i d(x_0,h_i(x_0))\]
We have for all $x\in X$
\[d(x_0,h_i(x))\leq M_0 + d(h_i(x_0),h_i(x))\]
By \eqref{eqn:coa} and Lipschitz property of $h_i$ we conclude that there exists $M_1>0$ such that for $d(x_0,x)>M_1$

\begin{equation} \label{eqn:r}
\sum p_i(x) \log^+(d(x_0,h_i(x))) \leq (1-b/2) \log(d(x_0,x))
\end{equation}
Let 
\[M_2=\sup_{d(x,x_0)<M_1} \max_i d(h_i(x),x_0)\]
We divide $\tilde{\nu}$ into two parts: $\tilde{\nu}_1$, supported on a ball $B_{M_1}(x_0)$ and $\tilde{\nu}_2$, supported on the outside.
Clearly, 
\[F(\tilde{\nu})=F(\tilde{\nu_1})+F(\tilde{\nu}_2)\]
We have
\[\int_X \log^+(d(x,x_0)) F(\tilde{\nu}_1)(dx) \leq \log(M_2)\]
and (by \eqref{eqn:r})
\[\int_X \log^+(d(x,x_0)) F(\tilde{\nu}_2)(dx) \leq (1-b/2) \int_X \log^+(d(x,x_0)) d\tilde{\nu}(x)\] 
Hence, 
\[\int_X \log^+(d(x,x_0)) F(\tilde{\nu})(dx) \leq \log(M_2)+ (1-b/2) \int_X \log^+(d(x,x_0)) d\tilde{\nu}(x)\]
and the upper limit of $\int_X \log^+(d(x,x_0)) F^n(\tilde{\nu})(dx)$ is not greater than $2\log(M_2)/b$.
At the same time, by \cite{BDEG} the sequence $F^n(\tilde{\nu})$ converges weakly to $\nu$.
\end{proof}
\end{lem}

Let $h$ be the supremum of $h_i'(x)$ over all $i$ and all $x$.
Let 
\[A=\{x_0,h_1(x_0),\ldots,h_k(x_0)\}\]
Choose two positive constants $b''<b'<b$ and let $\delta>1$ be such that
\begin{equation} \label{eqn:delta}
b'-(\delta -1) \log h > b'' \delta
\end{equation}

Let $n_0=\lceil (\delta-1)^{-1} \rceil$, $n_j=\lfloor \delta n_{j-1} \rfloor$ and 
\[Z_j=\{\omega \in \Sigma; \diam (h_{\omega^{-n_j}}(A)) < \exp(-n_j b')\}\]

Similarly, let
\[\tilde{Z}_j=\{\omega \in \Sigma; \diam (h_{\omega^{n_j}}(A)) < \exp(-n_j b')\}\]

Both $Z_j$ and $\tilde{Z}_j$ are disjoint unions of $n_j$-th level cylinders, hence by \eqref{eqn:mumi}
\[\mu_-(Z_j)=\mu_+(\tilde{Z}_j)\]
By Lemma \ref{lem:clt},
\[p_x(\tilde{Z}_j) \geq 1-\frac {L(A)(1+\log^+(d(x,x_0)))} {(b-b')^2 n_j}\]
Using Lemma \ref{lem:log} we get
\[\mu_+(\tilde{Z}_j)=\int_X p_x(\tilde{Z}_j) d\nu(x) \geq 1-\frac {L(A)(1+M(x_0))} {(b-b')^2 n_j}\]
As $n_j$ diverge exponentially fast, the series $\sum (1-\mu_-(Z_j))=\sum (1-\mu_+(\tilde{Z}_j))$ is convergent.
Denote $Z=\liminf Z_j$.
By Borel-Cantelli Lemma
\[\mu_-(Z)=1\]

As 
\[\diam (h_j\circ h_{\omega^n}(A))\leq h\diam(h_{\omega^n}(A)),\] 
\eqref{eqn:delta} implies that for $\omega \in \bigcap_{j\geq m} Z_j$ and $n>n_m$
\[\diam (h_{\omega^{-n}}(A)) < \exp(-n b'')\]
In particular, for all $\omega\in Z$ the sequence $\diam (h_{\omega^{-n}}(A))$ is summable.

We can estimate
\[d(h_{\omega^{-(n+1)}}(x_0),h_{\omega^{-n}}(x_0)) = d(h_{\omega^{-n}}(h_{\omega_{n+1}}(x_0)),h_{\omega^{-n}}(x_0)) \leq \diam(h_{\omega^{-n}}(A))\]
hence $h_{\omega^{-n}}(x_0)$ form a Cauchy sequence.
We proved thus that the limit in \eqref{eqn:pidef} exists for all $\omega\in Z$, i.e. for a set of full measure $\mu_-$.
As stated above, \eqref{eqn:mudef} follows and we are done with Proposition \ref{prop:pidef}. \qed

Consider the space $X\times \Sigma$ with measure 
\[d\bar{\nu}(x,\omega) = p_x(d\omega) d\nu(x)\]
This measure is invariant under the map
\[g(x,\omega) = (h_{\omega_1}(x),\sigma \omega)\]
By Elton \cite{E}, Lemma 1, this map is ergodic.
An immediate corollary of Elton's Lemma and Proposition \ref{prop:pidef} is the following

\begin{prop} \label{prop:erg}
Under assumptions of Proposition \ref{prop:pidef}, the map $(\tilde{\Sigma}, \tilde{\sigma}, \mu)$ is ergodic.
\begin{proof}
Denote the ergodic average for any $f\in L^1(\tilde{\Sigma},\mu)$ by $\bar{f}$.
To prove the assertion we need to prove that $\bar{f}$ exists and is equal to $\mu(f)$ $\mu$-almost everywhere.

Note first that $(X\times\Sigma, g, \bar{\nu})$ is a factor of $(\Sigma\times\Sigma, S, \mu)$ under the projection $(\pi, \id)$.
Hence, the assertion is true for functions of the form $f(\tau,\omega)=f(\omega)$, i.e. depending only on the future.
As the ergodic average doesn't depend on the finite number of initial summands, same is true for functions of the form $f(\tau,\omega)=f(\tau^n,\omega)$ for any $n$.
Those functions are dense in $L^1(\sigma\times\Sigma, \mu)$, hence the assertion follows for all functions by standard density argument.
\end{proof}
\end{prop}

We remind that as $(\tilde{\Sigma}, \tilde{\sigma})$ is bijective, the map $(\tilde{\Sigma},\tilde{\sigma}^{-1}, \mu)$ is also ergodic.

\section{Ergodic properties of typical orbits}

Let us define the Lyapunov exponent and (minus) entropy of the system:
\[\lambda = \int \sum_{j=1}^k p_j(x) \log|h_j'(x)| d\nu(x),\]
\[\eta = \int \sum_{j=1}^k p_j(x) \log p_j(x) d\nu(x)\]
By \eqref{eqn:coa}, $\lambda<0$.

For any $z_1,\ldots,z_k\in L^1(\R^d,\nu)$ we denote 
\[z=\int_{\R^d} \sum_{j=1}^k p_j(x) z_j(x) d\nu(x)\]
and 
\[z_{\omega^n}(x) = \sum_{i=1}^n z_{\omega_i}(h_{\omega^{i-1}}(x)).\]

The following two lemmas are direct consequences of the Birkhoff Ergodic Theorem; the first by application to the map $(\tilde{\Sigma}, \tilde{\sigma}, \mu)$ and the function $f(\tau,\omega)=z_{\omega_1}(\pi(\tau))$, the second is obtained for the same $f$ but for map $(\tilde{\Sigma},\tilde{\sigma}^{-1}, \mu)$.
Note that 
\[\mu(f)= \bar{\nu}(f\circ (\pi^{-1},\id)) = \int_{\R^d} \sum_{j=1}^k z_j(x) p_j(x) d\nu(x) =z\]

\begin{lem} \label{lem:est}
For any $\epsilon>0$, $\mu_+$-almost every $\omega\in\Sigma$ and $\nu$-almost all $x\in\R^d$ one can find a constant $K>0$ such that for all $n$
\[K^{-1} \exp (n(1+\epsilon)z)\leq \exp(z_{\omega^n}(x))\leq K \exp (n(1-\epsilon)z).\]
\end{lem}

\begin{lem} \label{lem:est2}
For any $\epsilon>0$ and $\mu_-$-almost every $\omega\in\Sigma$ there exists $K>0$ such that for all $n$
\[K^{-1} \exp (n(1+\epsilon)z)\leq \exp(z_{\omega^{-n}}\circ \pi(\sigma^n\omega))\leq K \exp (n(1-\epsilon)z).\]
\end{lem}

The following is a standard bounded distortion lemma.

\begin{lem} \label{lem:bdp}
For every $\epsilon>0$ and $K>0$ there exists $L(K)$ such that if for some $x\in \R^d$ and $\omega\in\Sigma$ and for all $n>0$
\begin{equation} \label{eqn:lem}
K^{-1} \exp (n(1+\epsilon)\lambda)\leq h_{\omega^n}'(x)\leq K \exp (n(1-\epsilon)\lambda)
\end{equation}
then
\[\sum_{m=0}^n \diam(h_{\omega^m}(B_{L(K)}(x))) \leq 1\]
\begin{proof}
Choose some $L$.
We have
\begin{equation} \label{eqn:zz}
\diam h_{\omega^n}(B_L(x)) \leq L \cdot \sup_{y\in B_L(x)} h_{\omega^n}'(y)
\end{equation}
As $\log |h_i'|$ is $C^\alpha$, for all $y\in B_L(x)$

\begin{equation} \label{eqn:z}
|h_{\omega^n}'(y)-h_{\omega^n}'(x)|\leq c \sum_{i=0}^{n-1} (d(h_{\omega^i}(x),h_{\omega^i}(y)))^\alpha h_{\omega^n}'(x)
\end{equation}
Denote 
\[d_0=2L\]
and
\[d_m=2LK \exp(m(1-\epsilon)\lambda) (1+c\sum_{j=0}^{m-1} d_j^\alpha)\]
By \eqref{eqn:zz} and \eqref{eqn:z} we have
\[d_m\geq \diam (h_{\omega^m}(B_L(x)))\]

Consider the series
\[M(L)=\sum_{i=0}^\infty (L(1+c)K \exp (i (1-\epsilon) \lambda))^\alpha\]
For $L$ small enough (depending on $K$) $M(L)$ is smaller than 1.
For such $L$ for all $n$
\[\sum_{m=0}^n d_m^\alpha\leq M(L)\]
It implies that
\[\sum_{m=0}^n \diam(h_{\omega^m}(B_{L(K)}(x))) \leq 2LK (1+cM(L)) \sum_{m=0}^n \exp(m(1-\epsilon)\lambda) \leq\] 
\[\leq \frac {2LK (1+cM(L))} {e^{-\lambda} -1}\]
and the right hand side is arbitrarily small (independent of $n$) for $L$ sufficiently small.
\end{proof}
\end{lem}

\begin{cor} \label{cor:bdp}
Under assumptions of Lemma \ref{lem:bdp}, for any $C^\alpha$ functions $z_i$ the difference $|z_{\omega^n}(x)-z_{\omega^n}(y)|$ is uniformly bounded for all $y\in B_{L(K)}(x)$.
\end{cor}

The following bounded distortion lemma implies that both \eqref{eqn:str}
and\eqref{eqn:reg} are preserved under iterations (up to rescaling).

\begin{lem} \label{lem:bdp2}
Fix $x\in \R^d$, $R>0$ and $\omega^n$.
Assume that for some $K>0$ and $\epsilon>0$ and for all $m\leq n$
\[K^{-1} \exp(m (1+\epsilon) \lambda) \leq h_{(\sigma^{n-m}\omega)^m}'(h_{\omega^m}(x)) \leq K \exp (m (1-\epsilon) \lambda)\]
Then there exists $l(K,R)$ such that for all $r\leq l(K) \exp(n(1+2\epsilon)\lambda)$ we have
\[\sum_{m=0}^n \diam(h_{(\sigma^{n-m}\omega)^m}^{-1}(B_r(h_{\omega^n}(x)))) \leq R\] 
\begin{proof}
Choose some $l$ and let 
\[r= l \exp(n(1+2\epsilon)\lambda)\]
Like in the proof of Lemma \ref{lem:bdp}, we can estimate 
\[\diam(h_{(\sigma^{n-m}\omega)^m}^{-1}(B_r(h_{\omega^n}(x)))) \leq d_m\]
where
\[d_0=2r\]
and
\[d_m= 2r K \exp(-m (1+\epsilon) \lambda)(1+c \sum_{j=0}^{m-1} d_j^\alpha)\]

Consider the series
\[M(l)=\sum_{j=0}^\infty (2 l K (1+c) \exp(j \epsilon \lambda))^\alpha\]
For $l$ small enough (depending on $K$), 
\[M(l)<1\]
and
\[\sum_{m=0}^{n-1} d_m^\alpha < M(l)\]
Hence, 
\[\sum_{m=0}^n \diam(h_{(\sigma^{n-m}\omega)^m}^{-1}(B_r(h_{\omega^n}(x))))\leq 2lK (1+M(l))e^{n (1+2\epsilon)\lambda} \sum_{m=0}^n e^{-m (1+\epsilon) \lambda} \leq\] 
\[\leq \frac {2e^{-\lambda} l K (1+c)} {e^{-(1+\epsilon)\lambda}-1}\]
and for $l$ small enough the right hand side is smaller than any constant, independently of $n$.
\end{proof}
\end{lem}

\begin{cor} \label{cor:bdp2}
Under assumptions of Lemma \ref{lem:bdp2}, for any $C^\alpha$ functions $z_i$ the difference $z_{\omega^n}(x)-z_{\omega^n}(y)$ is uniformly bounded for all $y\in h_{\omega^n}^{-1}(B_r(x))$.
\end{cor}

Before we begin the proof of Theorem \ref{thm:main}, let us show an example of iterated function system satisfying the regular open set condition.
The maps are not similitudes (it is impossible to have OSC for similitudes that are not all contracting) but their derivatives are constant in $U$.
Outside $U$ it may well happen that all the maps are expanding at some point, hence \eqref{eqn:coa} is satisfied only on $U$.
The system is thus not necessarily contracting-on-average.
However, Theorem \ref{thm:main} still holds.

\begin{ex}
Let $B_n=(10^n,3\cdot 10^n)$ for $n\geq 0$.
Let 
\[h_1(x)=\frac 1 {20} x + 15 \cdot 10^{n-2}\]
for $x\in B_n$, $n>0$,
\[h_1(x)=\frac 1 {20} x + 1\]
for $x\leq 3$,
\[h_2(x)= 5 x + 5 \cdot 10^n\]
for $x\in B_n$, $n>0$,
\[h_2(x)= 5 x + 5\]
for $x\leq 3$.
We extend maps $h_1,h_2$ to be homeomorphisms of $\R$.
The system satisfies the regular open set condition for $U=\bigcup B_n$ (regularity is satisfied because all the components of $U$ have diameters bounded away from 0).
The system is contracting-on-average (on $U$ only) for 
\[p_1(x)=p_1=1-p_2>\frac {\log(50/7)} {\log(500/17)}\]

For negative $x$ the contracting map $h_1$ is stronger than the dilatating map $h_2$ (in addition it is applied more often for typical $\omega$) and the fixed point for $h_2$ is smaller than the fixed point for $h_1$.
Similarly, for positive $x$ both $h_1$ and $h_2$ are equally strong (former moves $B_i$ into $B_{i-1}$, latter into $B_{i+1}$) but the former happens more often for typical $\omega$.
Hence, wherever we start, the typical trajectory lands inside $B_0$ after finite number of iterations.
The invariant measure (if it exists) is thus supported on $U$, hence independent from
the way we define $h_1$ and $h_2$ outside $U$.

Note that we may define $h_1$ and $h_2$ in such a way that the system
satisfies \eqref{eqn:coa} on whole $\R$. Hence, the invariant measure
indeed exists and is unique.
Its Hausdorff dimension equals
\begin{equation} \label{eqn:hd}
\dim_H(\nu)=\frac {p_1 \log p_1 +p_2 \log p_2} {-2p_1 \log 2 - (p_1-p_2)
\log 5}.
\end{equation}

More detailed study shows that the necessary and sufficient condition for existence of invariant measure is $p_1>1/2$, even though the iterated function system doesn't satisfy \eqref{eqn:coa} in this case.
The invariant measure is always unique and its Hausdorff dimension is still given by \eqref{eqn:hd}.
\end{ex}

\section{Proof of Theorem \ref{thm:main}}

\subsection{Upper bound}

Fix small $\epsilon>0$.
Choose some big $K>0$ and $n>0$.

Let $\Sigma(n,K,\epsilon)(x)$ be the set of sequences $\omega\in\Sigma$ for which 
\[K^{-1} \exp (n(1+\epsilon) \eta) \leq p_{\omega^n}(x) \leq K \exp (n(1-\epsilon)\eta)\]
and
\begin{equation} \label{eqn:3}
K^{-1} \exp (n(1+\epsilon)\lambda)\leq h_{\omega^n}'(x)\leq K \exp (n(1-\epsilon)\lambda).
\end{equation}

By Lemma \ref{lem:est}, $\mu(\Sigma(n,K,\epsilon)(x))$ is (for sufficiently big $K$, independently of $n$) greater than $(1-\epsilon)$ for all $x$ except some set $A(n,K,\epsilon)$ of $\nu$-measure smaller than $\epsilon$.
In addition, $\nu([-K,K]^d)>1-\epsilon$ for $K$ big enough.

Let us assume $K$ is sufficiently big that all the above is satisfied.
Note that (as $\sum_{\omega^n} p_{\omega^n}(x)=1$) $\omega^n$ can take
at most $K \exp(-n(1+\epsilon)\eta)$ values inside each $\Sigma(n,K,\epsilon)(x)$.
Given $\omega^n$, let $B(\omega^n,K,\epsilon)$ denote the set of $x$
for which $\omega^n=\tau^n$ for some $\tau\in \Sigma(n,K,\epsilon)(x)$.
Denote
\[\nu_n=\sum_{\omega^n} h_{\omega^n *}
(p_{\omega^n} \cdot \chi_{B(\omega^n,K,\epsilon) \cap [-K,K]^d}\cdot \nu)\]

As 
\[\nu = \sum h_{\omega^n *} (p_{\omega^n}\nu),\]
we have $\nu_n\leq \nu$.

At the same time, for any $x\in [-K,K]^d\setminus A(n,K,\epsilon)$ we have
\[p_x(\omega; x\in B(\omega^n, K, \epsilon)) \geq 1-\epsilon\]
Hence,
\[|\nu_n|\geq \nu([-K,K]^d\setminus A(n,K,\epsilon)) \cdot \inf_{x\in [-K,K]^d\setminus A(n,K,\epsilon)}p_x(\omega; x\in B(\omega^n,K,\epsilon)) \geq 1-3\epsilon.\]

Let $(E_i^{\epsilon,n})$ be a family of sets of diameter $L(K)$, covering $[-K,K]^d\setminus A(n,K,\epsilon)$.
This family may be chosen to have at most $(2dK/L(K))^d$ elements.
Let us choose a point $x_i^{\epsilon,n}$ inside each
$E_i^{\epsilon,n}\setminus A(n,K,\epsilon)$ whenever this set is nonempty.
Let $(F_j^{\epsilon,n})$ be a family of sets of the form
$h_{\omega^n}(E_i^{\epsilon,n})$ for $\omega\in
\Sigma(n,K,\epsilon)(x_i^{\epsilon,n})$.
This family has at most $2^d d^d K^{d+1} L(K)^{-d} \exp(-n(1+\epsilon)\eta)$ elements.
By Lemma \ref{lem:bdp}, Corollary \ref{cor:bdp} and \eqref{eqn:3} each of sets $F_j^{\epsilon,n}$ has diameter not greater than $K_0 L(K) K \exp (n(1-\epsilon)\lambda)$.

As the support of $\nu_n$ is contained in $\bigcup_{\omega^n} h_{\omega^n}([-K,K]^d\cap B(\omega^n,K,\epsilon))$, it is also contained in $\bigcup F_j^{\epsilon,n}$.
The sum
\[\sum (\diam F_j^{\epsilon,n})^{s(1+3\epsilon)}\leq 2^d d^d K_0^{s(1+3\epsilon)} K^{d+1+s(1+3\epsilon)} L(K)^{s(1+3\epsilon)-d} \exp(n \eta \epsilon(1-3\epsilon))\]
is arbitrarily small for big $n$.
Let us choose $n$ for which this sum is smaller than $\epsilon$ and denote $Z_\epsilon=\bigcup F_j^{\epsilon,n}$.

We can now repeat the procedure for a smaller $\epsilon_m=2^{-m}\epsilon$.
The set $Z=\limsup Z_{\epsilon_m}$ has full $\nu$-measure and we have constructed a family of its covers $(G^M)=\bigcup_{m>M} \bigcup F_j^{\epsilon_m}$ such that
\[\sum_{G_l^M\in G^M} (\diam G_l^M)^{s(1+3\epsilon)} \leq 1.\]
Hence, 
\[\dim_H(Z)\leq s (1+3\epsilon)\]
and $\epsilon$ can be taken arbitrarily small. \qed

\subsection{Lower bound}

We claim that if we can find a family of sets $\Sigma_K\in\Sigma$ such
that $\mu(\Sigma_K)\nearrow 1$ and the measures
$\nu_K=\pi_*(\mu_{|\Sigma_K})$ have lower Hausdorff dimension not
smaller than $s_0$ then the lower Hausdorff dimension of $\nu$ will be not smaller than $s_0$.

Suppose it is not true, i.e. there exists a set $X\subset\R^d$ such that $\nu(X)>0$ and $\dim_H(X)<s_0$.
For some $K$, $\nu(X)+\mu(\Sigma_K)>1$.
Hence, the set $Y=\pi^{-1}(X)\cap \Sigma_K$ has positive measure $\mu$.
It follows that $\pi(Y)$ has positive measure $\nu_K$ and its Hausdorff 
dimension is smaller than $s_0$ -- a contradiction.

Thus, we need only to prove the existence of such families of sets $\Sigma_K$ (or equivalently, of measures $\nu_K$) of lower Hausdorff dimension arbitrarily close to $s$.

Fix $\epsilon>0, K>0$.
Let $\Sigma_K$ be the set of $\omega\in \Sigma$ for which for all $n$
\begin{equation} \label{eqn:lb1}
K^{-1} \exp (-n(1-\epsilon)\lambda)\leq |(h_{\omega^{-n}}^{-1})'\circ \pi(\omega)|\leq K \exp (-n(1+\epsilon)\lambda)
\end{equation}
and
\begin{equation} \label{eqn:lb2}
K^{-1} \exp (n(1+\epsilon)\eta)\leq p_{\omega^{-n}}\circ \pi(\sigma^n\omega)\leq K \exp (n(1-\epsilon)\eta).
\end{equation}
For $K$ big enough $\mu_-(\Sigma_K)$ is arbitrarily close to 1 by Lemma \ref{lem:est2}.

Let $x=\pi(\omega)$ for some $\omega\in\Sigma_K$ and
\[r_n = K^{-1} l(K,R) \exp(n(1+2\epsilon)\lambda),\]
where $R=R_1$ in SOSC case or $R=R_2$ in ROSC case.
The ball $B_{r_n}(x)$ may in general intersect many of the sets $h_{\tau^n}(U)$.
Let $T$ be the set of all $\tau^n$ for which $\nu_K(B_{r_n}(x)\cap h_{\tau^n}(U))>0$.

By Lemma \ref{lem:bdp2}, Corollary \ref{cor:bdp2} and \eqref{eqn:lb2} we have
\[p_{\tau^n}(h_{\tau^n}^{-1}(y)) \leq K_0 K \exp(n(1-\epsilon)\eta)\]
for $\nu_K$-almost all $y\in B_{r_n}(x)\cap h_{\tau^n}(U)$.
Hence,
\[\nu_K(h_{\tau^n}(U)\cap B_{r_n}(x))\leq \nu(h_{\tau^n}(U)\cap B_{r_n}(x)) \leq K_0 K \exp (n (1-\epsilon) \eta)\]
for all $\tau^n\in T$.

We claim that $T$ has uniformly bounded (independently of $n$) number of elements.
Assume first that the strong open set condition is satisfied. 
Choose any point $y\in U$ and any $\tau^n\neq \omega^n$.
Let $m=\min \{j\geq 1; \omega_j\neq \tau_j\}-1$.
By SOSC, 
\[d(h_{(\sigma^m\omega)^{-(n-m)}}(h_{\omega^{-n}}^{-1}(x)), h_{(\sigma^m \tau)^{-(n-m)}}(y))>R_1\]
Hence, by Lemma \ref{lem:bdp2} 
\[h_{\tau^{-n}}(y)\notin B_{r_n}(x)\]
and thus $\tau^{-n} \notin T$.
Hence, $T$ has just one element $\omega^{-n}$.

Consider now the regular open set condition situation. 
Choose $\tau^n\in T$ and $\nu_K$-typical $y\in h_{\tau^n}(U)\cap B_{r_n}(x)$.
By Lemma \ref{lem:bdp2} and Corollary \ref{cor:bdp2}, 
\[B_{K_0r}(h_{\tau^n}^{-1}(y))\subset h_{\tau^n}^{-1}(B_{r_n}(x))\subset B_r(h_{\tau^n}^{-1}(y))\]
for some $r<R_2$.
By ROSC
\[\vol(U\cap h_{\tau^n}^{-1}(B_{r_n}(x))) \geq R_3 K_0^d \vol(h_{\tau^n}^{-1}(B_{r_n}(x)))\]
and (by Corollary \ref{cor:bdp2} again)
\[\vol(h_{\tau^n}(U) \cap B_{r_n}(x)) \geq K_1 \vol(B_{r_n}(x))\]
As $h_{\tau^n}(U)$ are pairwise disjoint, $T$ can have at most $K_1^{-1}$ elements.

Hence, if either SOSC or ROSC are satisfied,
\[\nu_K(B_{r_n}(x)) \leq K_2 \exp (n (1-\epsilon) \eta)\]
with $K_2$ depending on $K$ but not on $x$ or $n$.
Hence, 
\begin{equation} \label{eqn:eee}
\limsup_{r\rightarrow 0} \frac {\nu_K(B_r(x))} {(2r)^{s(1-4\epsilon)}} =0
\end{equation}
We only assumed above \eqref{eqn:lb1} and \eqref{eqn:lb2}, so \eqref{eqn:eee} is satisfied $\nu_K$-almost everywhere, which (by Frostman Lemma) implies that the lower Hausdorff dimension of $\nu_K$ is not smaller than $s(1-4\epsilon)$.\qed

\newpage
\bibliography{ref}

\end{document}